\documentclass[10pt, a4paper]{article}

\usepackage{color}
\usepackage{bbm}
\usepackage{amssymb,amsmath,amsthm}
\usepackage{verbatim}
\usepackage{setspace}
\usepackage{graphicx}
\usepackage{epsfig}
\normalsize \evensidemargin=0pt \oddsidemargin=0pt \hoffset=6pt
\newcommand\smx[1]{\left(\begin{smallmatrix}#1\end{smallmatrix}\right)}

\newtheorem{thm}{Theorem}[section]

\newtheorem*{rem}{Remark}

\usepackage{hyperref}

\begin{document}

\title{Upper and lower bounds for the reliability measure of a discrete distribution conditionally on the first three moments}

\author{DAVIDE DI CECCO\footnote{email: {\tt davide.dicecco@gmail.com}}}
\date{}
\maketitle

\begin{abstract}
We give sharp bounds for the reliability measure of a discrete
r.v. defined on $\{0,\ldots,n\}$, conditionally on the knowledge
of the first three moments of the r.v. The present work is as an
extension of the results given in [Di Cecco, Stat. Prob. Lett.,
\textbf{81}(2011), 411--416].
\end{abstract}

\section{Introduction}

Let $S$ be a discrete r.v. defined on the integers
$\{0,\ldots,n\}$ and $P$ be its distribution. Let $R_{k,n}$
denotes the probability $P(S\geq k)$, and let $\mu_i$ be the
$i$--th  moment of $S$, $\mu_i=E[S^i]$. In the present paper we
give exact upper and lower bounds for $R_{k,n}$, for any $k$ $\in$
$\{0,\ldots,n\}$, conditionally on the knowledge of the first
three moments of $S$, $\mu_1$, $\mu_2$, and $\mu_3$. Additionally,
we define the two extremal distributions on $\{0,\ldots,n\}$
consistent with the given $\mu_1$, $\mu_2$ and $\mu_3$, achieving
the maximum and the minimum for $R_{k,n}$.

A finite discrete distribution on $\{0,\ldots,n\}$ is completely
identified by its first $n$ moments. In fact, the distribution of
$S$ (and $R_{k,n}$ too), can be written in terms of its first $n$
factorial moments,
$\widetilde{\mu}_i=E\left[\binom{S}{i}i!\right]$, $i=1,\ldots,n$,
(see, e.g., \cite{Prekopa90}):
\begin{equation}\label{eq: JK}
P(S=k) = \sum_{i=k}^{n} (-1)^{i-k} \binom{i}{k}
\frac{\widetilde{\mu}_i}{i!},   \qquad   R_{k,n} = \sum_{i=k}^{n}
(-1)^{i-k} \binom{i-1}{k-1} \frac{\widetilde{\mu}_i}{i!},
\end{equation}
but the following relations hold between ordinary moments and
factorial moments:
\begin{equation}\label{eq: Stirling}
\mu_j = \sum_{i=1}^j S2(j,i) \, \widetilde{\mu}_i,
\qquad\qquad\qquad   \widetilde{\mu}_j = \sum_{i=1}^j S1(j,i) \,
\mu_i,
\end{equation}
where $S1(j,i)$ and  $S2(j,i)$  are the Stirling numbers of the
first and of the second kind respectively. Then, the distribution
of $S$ can be parameterized in terms of $\mu_1,\ldots,\mu_n$.

For any r.v. $S$ defined on $\{0,\ldots,n\}$ we have (or we can
construct) a sequence $(X_1,\ldots,X_n)$ of exchangeable Bernoulli
r.v.s, such that $S$ can be viewed as the tally variable of a set
of $n$ events $(S=\sum_{i=1}^n X_i)$. There is an obvious one to
one relation between their distribution which is given by
\[
P(S=k) = \binom{n}{k} P(X_1=x_1,\ldots,X_n=x_n) \qquad \text{with
} \sum_{i=1}^n x_i =k,
\]
and we can equivalently refer to the sequence or to the counting
variable.

Many ways to parameterize the joint distribution of $n$
exchangeable Bernoulli r.v.s $(X_1,\ldots,X_n)$ (and hence the
distribution of $S$), have been explored. De~Finetti since his
earlier works introduced a parameterization in terms of the
parameters $(w_1,\ldots, w_n)$ where
\[
w_i = P(X_{1}=1, \ldots  ,X_{i}=1) = E[X_{1} \cdots X_{i}].
\]
Bahadur in \cite{Bahadur61} introduced a parameterization in terms
of the generalized correlations $(\rho_2,\ldots,\rho_n)$ where
\[
\rho_{i}=\frac{E[(X_{1}-w_1)\cdots(X_{i}-w_1)]}{[w_1(1-w_1)]^{\frac{i}{2}}}\:.
\]
There is a one to one relation between the first $m$ elements of
each one of the three parameterizations: $(w_1,\ldots,w_m)$,
$(w_1,\rho_2,\ldots,\rho_m)$, and $(\mu_1,\ldots,\mu_m)$ for any
$m \in \{1,\ldots,n\}$. In fact, we have that
$\widetilde{\mu}_i=\binom{n}{i}i! w_i$, then, by \eqref{eq:
Stirling}, the relations between the parametrs $w_i$ and the parameters $\mu_i$ are immediately derived; while the relations between the parameters $w_i$
and the parameters $\rho_i$ can be found in \cite{DiCecco11}

 As a consequence of these
relations, we have the equivalence of the bounds of $R_{k,n}$
conditioned on the knowledge of the first $m$ parameters of any of
the three parameterizations:
\[
\frac{\max}{\min}\:R_{k,n}(w_1,\ldots,w_m) =
\frac{\max}{\min}\:R_{k,n}(w_1,\rho_2,\ldots,\rho_m) =
\frac{\max}{\min}\:R_{k,n}(\mu_1,\ldots,\mu_m).
\]

In \cite{Zaigraev10} sharp bounds for $R_{k,n}(w_1)$ are given,
and in \cite{DiCecco11} sharp bounds for $R_{k,n}(w_1,\rho_2)$ are
given. In the present paper, in order to extend those results to
the first 3 parameters, we use the moment parameterization,
essentially for two reasons: firstly, formulae for the bounds of
$R_{k,n}$ in terms of $(\mu_1,\mu_2,\mu_3)$ have been proved to be
simpler; secondly, we want to show a link with the existing
literature that mostly refers to the moment parameterization,
rather than the other two parameterizations.

The result we obtain has an immediate interpretation in
Reliability theory, as $R_{k,n}$ represents the reliability of a
$k$--out--of--$n$ system, i.e., the probability that, in a system
of $n$ exchangeable components, at least $k$ will function.
Approximated bounds for the reliability measure of a discrete
distribution conditioned on its first (binomial) moments are
obtained via linear programming in  \cite{Prekopa90}. For an
analogous result on the reliability function of a continuous
distribution conditioned on the first moments see
\cite{CourtoisDenuit07}. Another example of a possible application
of the presented result involves developments of Condorcet's Jury
Theorem (see \cite{kaniovski09} and references therein), studying
the scenario of dichotomous voting in a jury (group of experts)
with a majority voting rule and certain hypotheses of dependence
among the jurors.

The present paper is a direct extension of the geometric approach
described in \cite{DiCecco11} (which, by the way, can be easily
employed to find sharp bounds for the probability $P(S=k)$ of $S$
being \textit{exactly} equal to $k$), but it can be read on its
own. In Section~\ref{sec: geometry} we present that geometric
approach, in Section~\ref{sec: result} we use it to state our
result.

\section{Some Geometry}\label{sec: geometry}

To outline our geometric approach, it is convenient to introduce
some notation. The convex hull of a set of points will be denoted
in angle brackets: $\langle\cdot \rangle$. Let $y_1,\ldots,y_m$ be
points in a $(m-1)$--dimensional space, where
$y_i=(y_{i,1},\ldots,y_{i,m-1})^{T}$; then $||y_1,\ldots,y_m||$
will denote the following:
\[
||y_1,\ldots,y_m|| = \det\left(%
\begin{array}{cccc}
   1   & y_{1,1} & \cdots & y_{1,m-1} \\
\vdots &  \vdots &        &\vdots      \\
   1   & y_{m,1} & \cdots & y_{m,m-1} \\
\end{array}%
\right).
\]
Let $p$ be a point in the $(m-1)$--dimensional space: when $p$ is
variable, $||y_1,\ldots,y_{m-1},p||=0$ is the equation of the
hyperplane $H$ containing points $y_1,\ldots,y_{m-1}$; while, if
$p$ is a fixed point, the sign of the determinant
$||y_1,\ldots,y_{m-1},p||$ reveals in which side of $H$ $p$
lies.
\begin{rem} In order to simplify the formulae that we are going to
describe, we normalize $S$, dividing it by $n$. So, from now on,
$\mu_j$ will denote the $j$-th moment of the r.v. $S/n$: $\mu_j =
E\left[\left(S/n\right)^j\right]$. In this manner, $\mu_j \in
[0,1]$, $\forall j$, and $\mu_1=w_1$.
\end{rem}
In order to find exact bounds for $R_{k,n}$ given
$(\mu_1,\ldots,\mu_m)$, we will consider the space
$\Phi_{k,n}^{(m)}$ of the admissible values for the array of
parameters $(\mu_1,\ldots,\mu_m,R_{k,n})$. We will show that
$\Phi_{k,n}^{(m)}$ is a bounded convex polytope of $m+1$ affine
dimensions, and hence we can calculate the maximum and the minimum
of $R_{k,n}$ given $(\mu_1,\ldots,\mu_{m})$ by finding the two
points of intersection of the vertical line $L_{\mu}$ passing
through $(\mu_1,\ldots,\mu_{m},0)$ and $(\mu_1,\ldots,\mu_{m},1)$,
with the upper and the lower boundaries of $\Phi_{k,n}^{(m)}$.

The space of the parameters $(\mu_1,\ldots,\mu_n)$, denote it
$\mathcal{M}_n$, is known to be an $n$--dimensional convex
polytope defined as the convex hull of the vertices
$\{v_{i,n}\}_{i=0,\ldots,n}$ (see, e.g., \cite{karlinStudden66}),
where
\[
v_{i,n} =
\left(\frac{i}{n},\left(\frac{i}{n}\right)^2,\ldots,\left(\frac{i}{n}\right)^n\right)^T.
\]
Denote as $\mathcal{M}_n^{(m)}$ the orthogonal projection of $\mathcal{M}_n$ over the first
$m$ axes. Obviously, $\mathcal{M}_n^{(m)}
= \langle v_{1,n}^{(m)},\ldots,v_{n,n}^{(m)}\rangle$ where
\[
v_{i,n}^{(m)} =
\left(\frac{i}{n},\left(\frac{i}{n}\right)^2,\ldots,\left(\frac{i}{n}\right)^m\right)^T.
\]
Let $c_d(t)$ be $(t,t^2,\ldots,t^d)^T$; the $d$--th order moment
curve is the curve parametrically defined as: $\{c_d(t) \:|\:
0\leq t\leq 1\}$. Both $\mathcal{M}_n$ and $\mathcal{M}_n^{(m)}$
are convex hulls of a set of points on the moment curve, hence are
cyclic polytopes (see, e.g., \cite{Grunbaum03}).

Each vertex $v_{i,n}$ of $\mathcal{M}_n$ represents the
distribution $S^*_i$ having $P(S^*_i=i)=1$. Under $S^*_i$, we have
$R_{k,n}=1$ if $i\geq k$, and $R_{k,n}=0$ if $i< k$. By \eqref{eq:
JK} and \eqref{eq: Stirling}, we can see that $R_{k,n}$ is a
linear function of $(\mu_1,\ldots,\mu_n)$. Define the point
$r_{i,k,n}^{(m)} \in \mathbb{R}^{m+1}$ as $(v_{i,n}^{(m)},0)$ if
$i<k$ and $(v_{i,n}^{(m)},1)$ if $i\geq k$; then, $\langle
r_{0,k,n}^{(m)},\ldots,r_{n,k,n}^{(m)}\rangle$ is exactly our
space $\Phi_{k,n}^{(m)}$.  So $\Phi_{k,n}^{(m)}$ is defined as the
convex hull of two sets of points lying on two parallel
hyperplanes (identified by $R_{k,n}=0$ and $R_{k,n}=1$): this kind
of convex polytope is sometimes called prismoid or prismatoid.
$\langle r_{0,k,n}^{(m)},\ldots,r_{k-1,k,n}^{(m)}\rangle$ is the
lower base of the prismatoid, call it $B_L$; $\langle
r_{k,k,n}^{(m)},\ldots,r_{n,k,n}^{(m)}\rangle=B_U$ is the upper base. The following theorem, whose proof is in Appendix, shed
some light on the structure of $\Phi_{k,n}^{(m)}$.
\begin{thm}\label{thm: m+2}
Any $m+2$ vertices  of $\Phi_{k,n}^{(m)}$ are affinely independent
unless they all belong to the same base.
\end{thm}
By Theorem \ref{thm: m+2}, each facet ($m$--dimensional face) of
$\Phi_{k,n}^{(m)}$, other than $B_L$ and $B_U$, is a simplex of
$m+1$ vertices. The projection of $\Phi_{k,n}^{(m)}$ over the
plane of the first $m$ axes is $\mathcal{M}_n^{(m)}$. In
particular, the projection of each facet of $\Phi_{k,n}^{(m)}$ is
a simplex inside $\mathcal{M}_n^{(m)}$, so, the projections of the
upper and the lower hulls of $\Phi_{k,n}^{(m)}$ provide two
subdivisions of $\mathcal{M}_n^{(m)}$. We call them the upper and
the lower subdivisions. Obviously, the given $(\mu_1,$ $\ldots,$
$\mu_{m})$ is a point of $\mathcal{M}_n^{(m)}$. To determine the
intersections of $L_{\mu}$ with the boundary of
$\Phi_{k,n}^{(m)}$, it suffices to find the two simplexes of the
upper and the lower subdivisions of $\mathcal{M}_n^{(m)}$
containing $(\mu_1,$ $\ldots,$ $\mu_{m})$, and then find the
intersections of $L_{\mu}$ with the two supporting hyperplanes of
$\Phi_{k,n}^{(m)}$ relative to the facets corresponding to those
simplexes. In the following Section \ref{sec: location}, we state
how to determine the simplexes containing $(\mu_1,\mu_2,\mu_3)$ in
the upper and lower subdivisions of $\mathcal{M}_n^{(3)}$.

\section{The main result}\label{sec: result}

\subsection{Point location}\label{sec: location}

For ease of notation, in the following we will denote
$r_{i,k,n}^{(3)}$ simply as $r_i$, and $v_{i,n}^{(3)}$ as $v_i$.
The following theorem, whose proof is in appendix, defines the
facial structure of $\Phi_{k,n}^{(3)}$.
\begin{thm}\label{thm: facets}
The upper facets of $\Phi_{k,n}^{(3)}$ are:
\[
\left\{%
\begin{array}{l}
\big\{\langle r_0,r_i,r_{i+1},r_k \rangle\big\}_{i=1,\ldots,k-2}; \\
\big\{\langle r_0,r_k,r_i,r_{i+1} \rangle\big\}_{i=k+1,\ldots,n-1}; \\
\big\{\langle r_t,r_{t+1},r_k,r_n \rangle\big\}_{t=0,\ldots,k-2}; \\
B_U =\big\{\langle r_k,r_t,r_{t+1}, r_n \rangle\big\}_{t=k+1,\ldots,n-2}. \\
\end{array}%
\right.
\]
The lower facets of $\Phi_{k,n}^{(3)}$ are:
\[
\left\{%
\begin{array}{l}
B_L =\big\{\langle r_0,r_i,r_{i+1}, r_{k-1} \rangle\big\}_{i=1,\ldots,k-3};\\
\big\{\langle r_0,r_{k-1},r_i,r_{i+1} \rangle\big\}_{i=k,\ldots,n-1};\\
\big\{\langle r_t,r_{t+1},r_{k-1},r_n \rangle\big\}_{t=0,\ldots,k-3};\\
\big\{\langle r_{k-1},r_t,r_{t+1},r_n \rangle\big\}_{t=k,\ldots,n-2}.\\
\end{array}%
\right.
\]
\end{thm}
$\mathcal{M}_n^{(3)}$ is a cyclic polytope, and its facial
structure is well known: its faces are $\big\{\langle
v_0,v_i,v_{i+1}\rangle\big\}_{i=1,\ldots,n-1}$ and $\big\{\langle
v_t,v_{t+1},v_n\rangle\big\}_{t=0,\ldots,n-2}$. We note in passing
that all the facets of $\Phi_{k,n}^{(3)}$ are of the kind $\langle
F,r_k\rangle$ or $\langle F,r_{k-1}\rangle$, for $F$ being a face
of $\mathcal{M}_n^{(3)}$.

Let $\mu$ be the point representing the array of the given
parameters $(\mu_1,\mu_2,\mu_3)$. The projections of the upper
facets of $\Phi_{k,n}^{(3)}$ form the upper subdivision of
$\mathcal{M}_n^{(3)}$, and can be divided into 4 groups of
simplexes which we will call blocks:
\[
\left\{%
\begin{array}{ll}
\big\{\langle v_0,v_i,v_{i+1},v_k \rangle\big\}_{i=1,\ldots,k-2}                &\text{block \textbf{1}};\\
\big\{\langle v_0,v_k,v_i,v_{i+1} \rangle\big\}_{i=k+1,\ldots,n-1}              &\text{block \textbf{2}};\\
\big\{\langle v_t,v_{t+1},v_k,v_n \rangle\big\}_{t=0,\ldots,k-2}                &\text{block \textbf{3}};\\
B_U^{(3)} =\big\{\langle v_k,v_t,v_{t+1}, v_n \rangle\big\}_{t=k+1,\ldots,n-2}  &\text{block \textbf{4}}.\\
\end{array}%
\right.
\]

All the simplexes in blocks \textbf{1} and \textbf{2} have the
edge $\langle v_0,v_k \rangle$ in common, so, if the point $\mu$
is in block \textbf{1} or \textbf{2}, to determine the simplex
$\langle v_0,v_k,v_{i^*},v_{i^* +1}\rangle$ containing it, we can
consider the dihedral angle $\xi$ between the two planes having
equations $||v_0,v_k,v_n,p||=0$ and $||v_0,v_k,v_{i^*},p||=0$. We
can calculate the cosine of $\xi$ and, equalling it to the cosine
of the dihedral angle between $||v_0,v_k,v_n,p||=0$ and
$||v_0,v_k,\mu,p||=0$ and solving for $i^*$, we obtain:
\begin{equation}\label{eq: i_blocks12}
i^* = \left\lfloor \frac{n(n \mu_3 - k \mu_2)}{n\mu_2-k\mu_1}
\right\rfloor ,
\end{equation}
where $\lfloor \cdot \rfloor$ denotes the floor function. All the
simplexes in blocks \textbf{3} and $B_U^{(3)}$ have the edge
$\langle v_k,v_n \rangle$ in common. Then, if $\mu$ is in block
\textbf{3} or $B_U^{(3)}$, we find the simplex $\langle
v_k,v_{t^*},v_{t^* +1},v_n\rangle$ containing it by equalling the
cosine of the dihedral angle  determined by the two planes
$||v_0,v_k,v_n,p||=0$ and $||v_{t^*},v_k,v_n,p||=0$ and the cosine
of the dihedral angle between $||v_0,v_k,v_n,p||=0$ and
$||v_{\mu},v_k,v_n,p||=0$. Solving for $t^*$ we obtain:
\begin{equation}\label{eq: t_blocks3U}
t^* = \left\lfloor \frac{n(n\mu_3-(k+n)\mu_2+
k\mu_1)}{n\mu_2-(k+n)\mu_1 + k} \right\rfloor.
\end{equation}

\subsection{The extremal distributions}

At this point, with few algebra,  we obtain the extremal
distribution $S^+$ on $\{0,\ldots,n\}$ consistent with $\mu_1$,
$\mu_2$ and $\mu_3$, achieving the maximum for $R_{k,n}$: we get
$\mu$ as a convex combination of the vertices of the simplex of
$\mathcal{M}_n^{(3)}$ containing it, and the coefficients of that
combination define $S^+$. The extremal distribution is clearly
unique, as none of the facets of $\Phi_{k,n}^{(3)}$ is orthogonal
to the plane of the first $3$ axes, so $L_{\mu}$ intersects the
upper boundary of $\Phi_{k,n}^{(3)}$ in a single point. $S^+$
concentrates the mass on four points, and, if $\mu$ is contained
in $\langle v_0,v_k,v_{i^*},v_{i^*+1}\rangle$, is defined as
\begin{equation}\label{eq: extremal1}
S^+=\left\{
\begin{array}{ll}
  0     & \text{with prob. }  \pi_{0}     = 1-\pi_k-\pi_{i^*}-\pi_{i^*+1};\\
  k     & \text{with prob. }  \pi_{k}     = \displaystyle\frac{ n[n^2\mu_3 - (2i^*+1)n\mu_2  + i^*(i^*+1)\mu_1 ]}{k(k-i^*)(k-i^*-1)};\\
  i^*   & \text{with prob. }  \pi_{i^*}   = \displaystyle\frac{ n[n^2\mu_3 - (k+i^*+1)n\mu_2 + k(i^*+1)\mu_1   ]}{i^*(k-i^*)};\\
  i^*+1 & \text{with prob. }  \pi_{i^*+1} = \displaystyle\frac{-n[n^2\mu_3 - (k + i^*)n\mu_2 + k i^* \mu_1     ]}{(i^*+1)(k-i^*-1)}.\\
\end{array}
\right.
\end{equation}
Moreover, if $\mu$ is in block \textbf{1} $(i^*<k)$, we have that
$\max(R_{k,n})=\pi_k$; if $\mu$ is in block \textbf{2} $(i^*
> k)$, $\max(R_{k,n})=\pi_k + \pi_{i^*}+ \pi_{i^*+1}$. When
$\mu$ is in $\langle v_k,v_{t^*},v_{t^*+1},v_n\rangle$, $S^+$ is
defined as
\begin{equation}\label{eq: extremal2}
S^+=\left\{
\begin{array}{ll}
  k     & \text{with prob. }  p_{k}     =\displaystyle\frac{ n[ n^2\mu_3 - n(n+2t^*+1)\mu_2  + ({t^*}^2+2nt^*+n+t^*)\mu_1 - (t^*+1)t^* ]}{(k-n)(k-t^*)(k-t^*-1)};\\
  t^*   & \text{with prob. }  p_{t^*}   =\displaystyle\frac{-n[ n^2\mu_3 - n(n+k+t^*+1)\mu_2 + (k+n+kn+kt^*+nt^*)\mu_1    - k-kt^*     ]}{(k-t^*)(n-t^*)};\\
  t^*+1 & \text{with prob. }  p_{t^*+1} =\displaystyle\frac{ n[ n^2\mu_3 - n(k+n+t^*)\mu_2   + (nt^*+kn+kt^*)\mu_1        - kt^*       ]}{(k-t^*-1)(n-t^*-1)};\\
  n     & \text{with prob. }  p_{n}     =1-p_k-p_{t^*}-p_{t^*+1}.\\
\end{array}
\right.
\end{equation}
in which case, if $\mu$ is in block \textbf{3} $(t^*<k)$, we have
that $\max(R_{k,n})=p_k + p_n$ and, if $\mu$ is in $B_U^{(3)}$
$(t^* > k)$, $\max(R_{k,n})=1$.

If $\mu$ is in block \textbf{1} or \textbf{2}, $i^*$, as defined
by \eqref{eq: i_blocks12}, is well defined, while $t^*$, as
defined by \eqref{eq: t_blocks3U}, can be out of the range
$\{0,\ldots,n\}$ and viceversa if $\mu$ is in block \textbf{3} or
in $B_U^{(3)}$. That is, one between $i^*$ and $t^*$ (but not
both) can be inadmissible, in which case, we can immediately state
that the simplex containing $\mu$ is the one determined by the
remaining value which is admissible. Otherwise, to determine which
of the two simplexes contains $\mu$, we can simply check which between
$(\pi_0,\pi_k,\pi_{i^*},\pi_{i^*+1})$ and
$(p_k,p_{t^*},p_{t^*+1},p_n)$ is a proper distribution. In fact,
one and one only of the two would be a set of values in $[0,1]$
summing to 1.

As regards to $\min(R_{k,n}(\mu_1,\mu_2,\mu_3))$, we proceed
similarly, by dividing the lower subdivision of
$\mathcal{M}_n^{(3)}$ into 4 groups of simplexes (blocks):
\[
\left\{%
\begin{array}{ll}
B_L^{(3)} =\big\{\langle v_0,v_i,v_{i+1},v_{k-1} \rangle\big\}_{i=1,\ldots,k-3}  &\text{block \textbf{1}};\\
\big\{\langle v_0,v_{k-1},v_i,v_{i+1} \rangle\big\}_{i=k,\ldots,n-1}           &\text{block \textbf{2}};\\
\big\{\langle v_t,v_{t+1},v_{k-1},v_n \rangle\big\}_{t=0,\ldots,k-3}             &\text{block \textbf{3}};\\
\big\{\langle v_{k-1},v_t,v_{t+1}, v_n \rangle\big\}_{t=k,\ldots,n-2}          &\text{block \textbf{4}}.\\
\end{array}%
\right.
\]
Then, the passages are the same as those of the upper bound, so we
limit ourselves to say that formulas \eqref{eq: i_blocks12} and
\eqref{eq: t_blocks3U}, are valid with $k-1$ substituting $k$, and
the corresponding extremal distribution $S^-$ is defined by
\eqref{eq: extremal1} and \eqref{eq: extremal2} with $k-1$
substituting $k$. In this case, if $\mu$ is in block \textbf{1}
$(i^*<k)$, we have that $\min(R_{k,n})=0$; if $\mu$ is in block
\textbf{2} $(i^*
\geq k)$, $\min(R_{k,n})=\pi_{i^*}+ \pi_{i^*+1}$; if $\mu$ is in block \textbf{3} $(t^*<k)$, $\min(R_{k,n})=p_n$; if $\mu$ is in block \textbf{4} $(t^* \geq k)$, $\min(R_{k,n})=p_{t^*}+ p_{t^*+1}+p_n$.

\begin{figure}[ht!]
\begin{center}
\begin{tabular}{ccc}
\epsfig{file=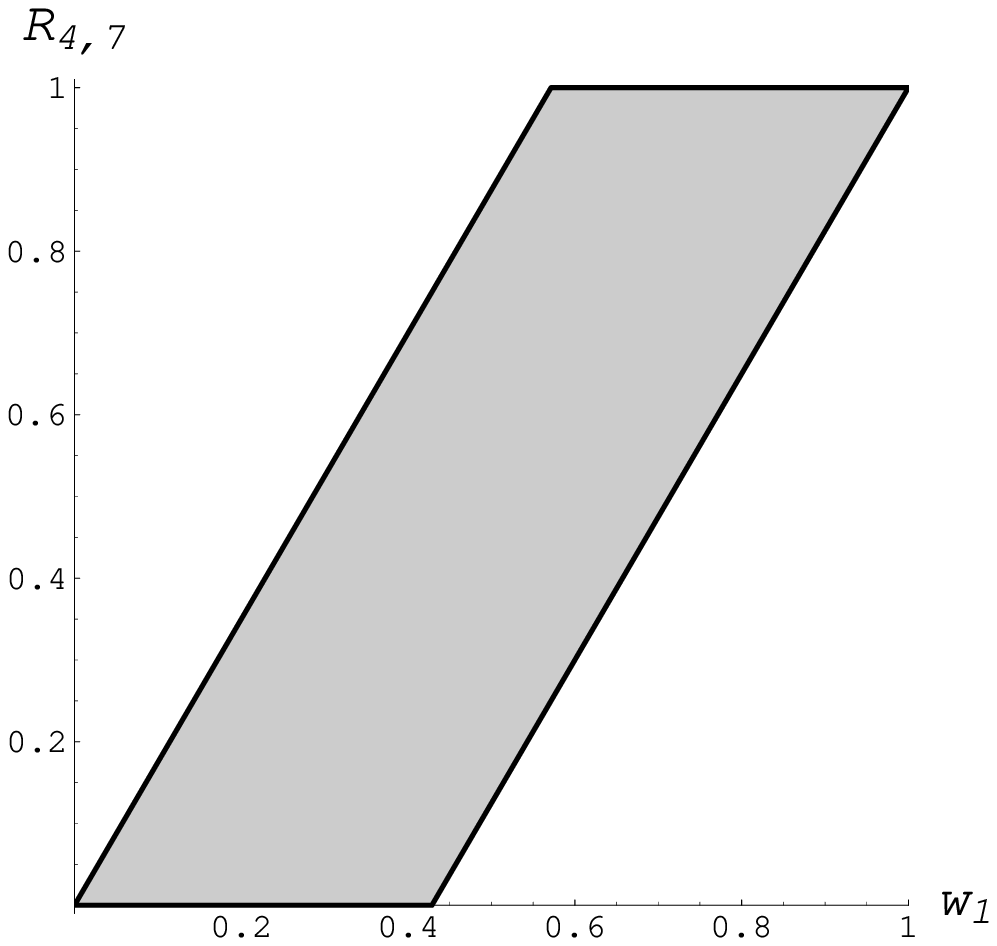, width=5cm}& \epsfig{file=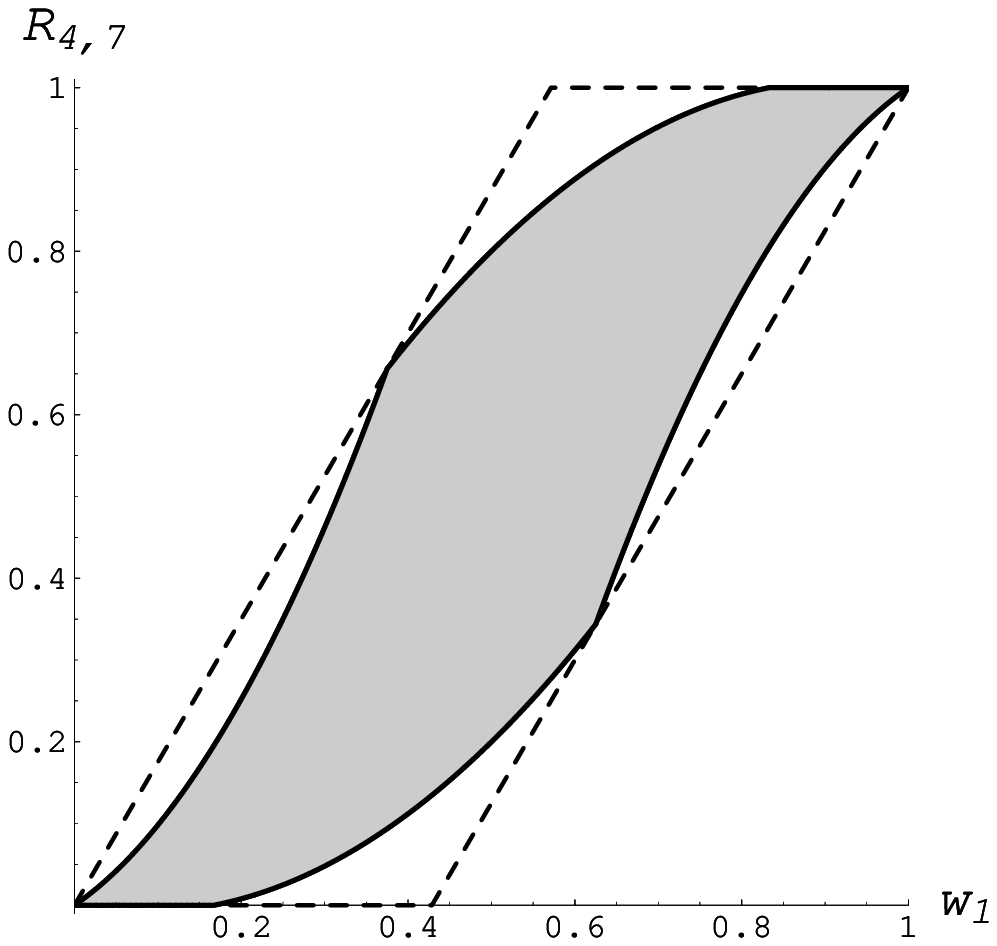,
width=5cm}& \epsfig{file=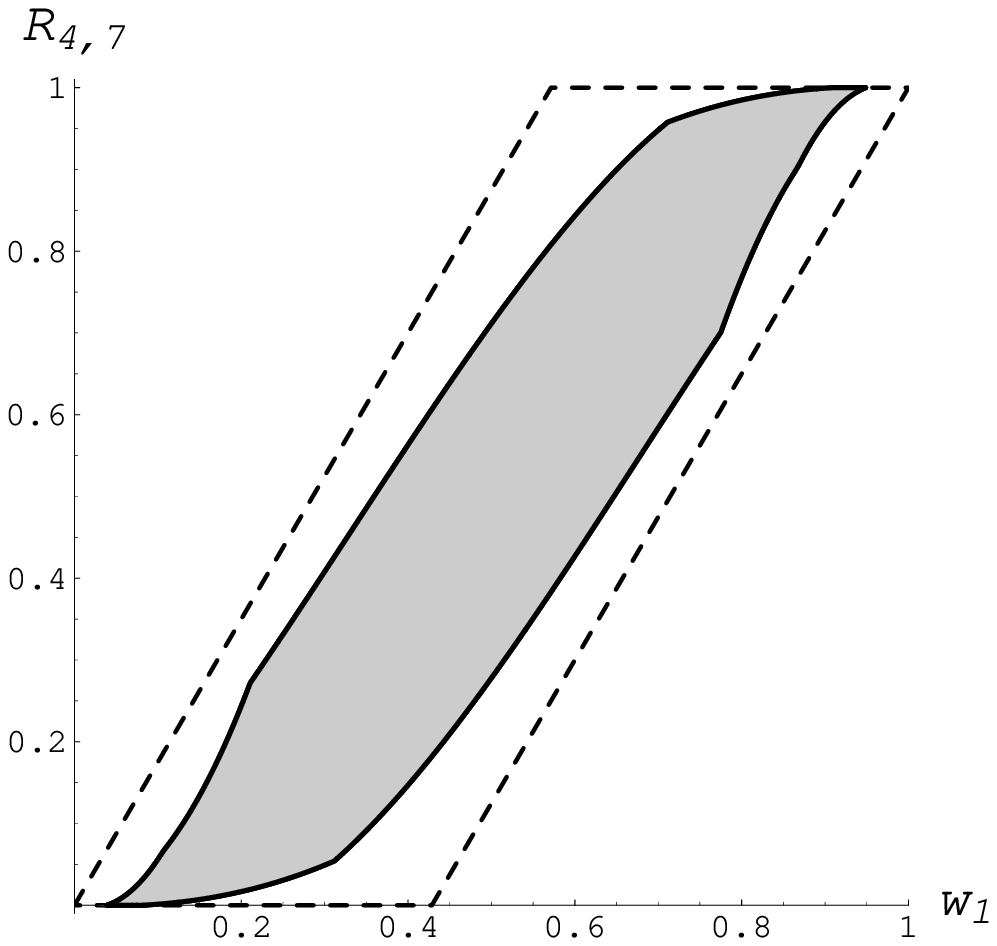, width=5cm}
\end{tabular}
\caption{Shaded areas represent the space of the admissible values
for $(w_1,R_{4,7})$ when no parameter is fixed (left), when we fix
$\rho_2=0.2$ (middle), and when we fix $\rho_2=0.2$ and
$\rho_3=0.1$ (right). Note that in the last case $\sim 0.03883
\leq w_1\leq \sim 0.94867 $.} \label{fig: plots}
\end{center}
\end{figure}

To give an example of the results obtained, we show in Figure
\ref{fig: plots} the bounds of $R_{k,n}$ as a function of $w_1$,
having fixed $\rho_2$ and $\rho_3$ (we condition on the
correlation parameters as they probably have a more interesting
interpretability than the moments). Note that, when we fix
$\rho_3$, $w_1$ cannot range freely in $[0,1]$, but has a narrower
interval of range which also depends on $w_1$, $\rho_2$, and $n$.

\section*{Appendix A: Proof of Theorem \ref{thm: m+2}}

\begin{figure}[ht!]
\begin{center}
\begin{tabular}{ccc}
\epsfig{file=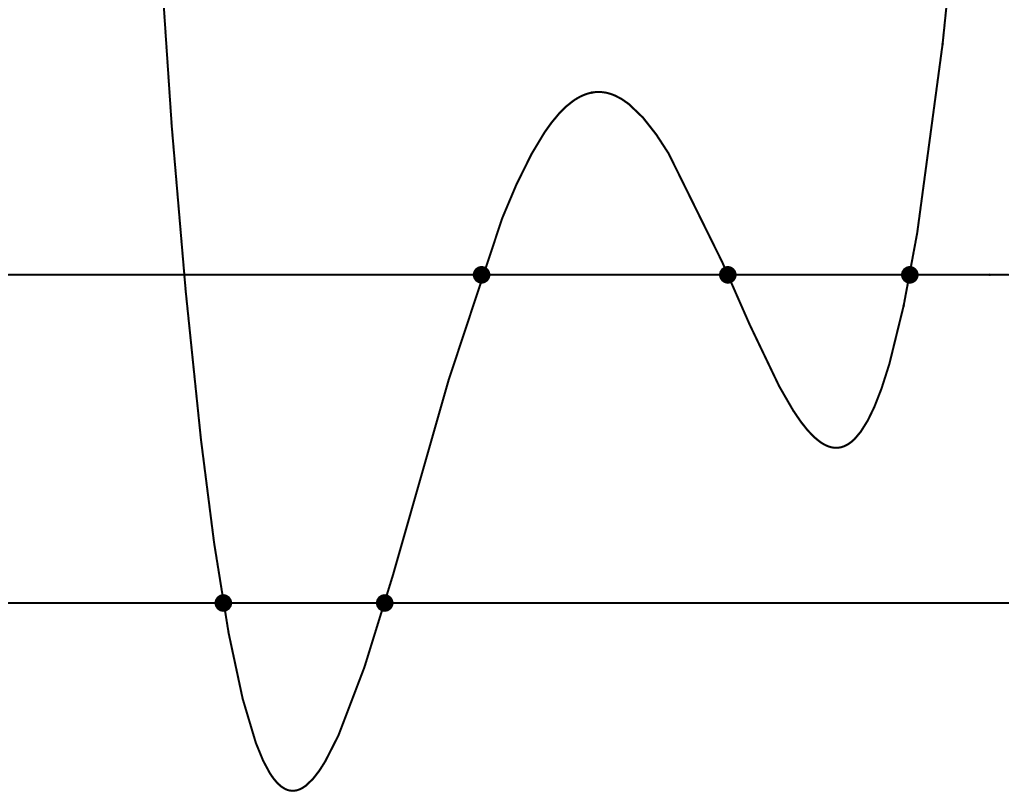, width=4cm} & \qquad\qquad\qquad &
\epsfig{file=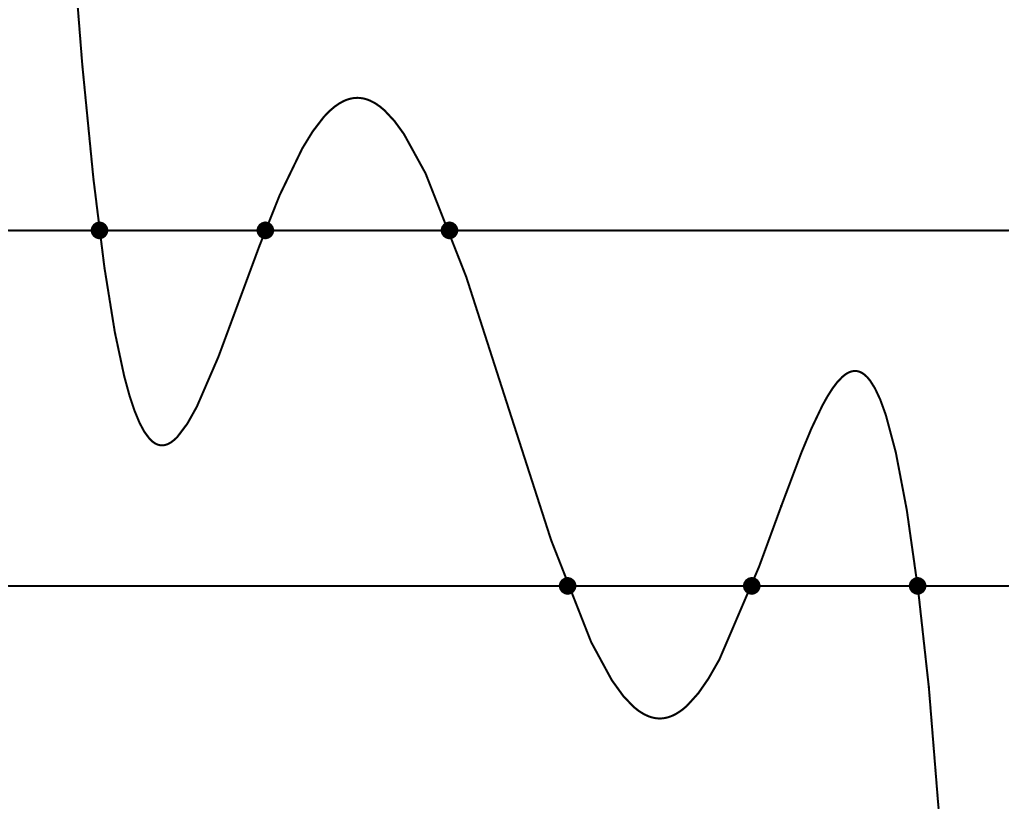, width=4cm}
\end{tabular}
\caption{An example of the disposition of the solution set when
$m=4$ (left) and $m=5$ (right).} \label{fig: proof}
\end{center}
\end{figure}

Suppose Theorem \ref{thm: m+2} is false. Then, $m+2$ vertices
$r_{i_1},\ldots, r_{i_{m+2}}$ of $\Phi_{k,n}^{(m)}$, which do not
belong all to the same base, lie on a common $m$--dimensional
hyperplane $H$. Fix, w.l.o.g., $i_1 < \ldots <i_t < k \leq i_{t+1}
< \ldots < i_{m+2}$. That is, $r_{i_1},\ldots, r_{i_{t}}$ belong
to $H \cap B_L$, while $r_{i_{t+1}},\ldots, r_{i_{m+2}}$ belong to
$H \cap B_U$. So, the two sets of vertices lie on two
$(m-1)$--dimensional parallel hyperplanes, and the same is valid
for their orthogonal projections over the first $m$ axes
(orthogonal projections maintain parallelism). That is, two sets
of vertices, $v_{i_1}, \ldots ,v_{i_{t}}$ and $v_{i_{t+1}}, \ldots
,v_{i_{m+2}}$, of $\mathcal{M}_n^{(m)}$ lie on two parallel
hyperplanes. Now, since the points $v_i$ belong to the moment
curve, there exist coefficients $a_0,\ldots,a_m$ such that the
equation
\begin{equation}\label{eq: H1}
a_0 + a_1 x + a_2 x^2 + \ldots + a_m x^m = 0
\end{equation}
has (at least) $t$ real roots $i_1,\ldots, i_t$;  and there exists
$b_0$, $b_0 \neq a_0$, such that the equation
\begin{equation}\label{eq: H2}
b_0 + a_1 x + a_2 x^2 + \ldots + a_m x^m = 0
\end{equation}
has at least $m+2-t$ roots $i_{t+1},\ldots, i_{m+2}$. Equations
\eqref{eq: H1} and \eqref{eq: H2} have at most $m$ real roots and
the corresponding polynomial curves have at most $m-1$ local
maxima and minima. It is easy to see  that, whether $m$ is odd or
even (Figure \ref{fig: proof} on the right and left respectively),
\eqref{eq: H2} can have (at most) $m+1-t$ roots greater than
$i_t$. So, $m+2$ vertices lying on two parallel hyperplanes cannot
exist.

\section*{Appendix B: Proof of Theorem \ref{thm: facets}}

By Theorem \ref{thm: m+2}, any 5 vertices of $\Phi_{k,n}^{(3)}$
are  linearly independent, hence, each facet other than $B_L$ and
$B_U$ is a simplex of 4 vertices, say $r_a$, $r_b$, $r_c$, $r_d$.
To find the facets, we search for the quadruples $(a,b,c,d)$ such
that the hyperplane $H$ passing through $r_a$, $r_b$, $r_c$, $r_d$
is a supporting hyperplane of $\Phi_{k,n}^{(3)}$, i.e., such that
all other vertices $r_x$ of $\Phi_{k,n}^{(3)}$ are on the same
side of $H$. That is, we search for the quadruples $(a,b,c,d)$
such that the determinant $||r_a,r_b,r_c,r_d,r_x||$, as a function
of $x$, has the same sign for all $x \in
\{0,\ldots,n\}\backslash\{a,b,c,d\}$. Fix $0\leq a<b<c<d\leq n$.
Three cases are possible:  
\[
\left\{%
\begin{array}{ll}
\textbf{1)} & a<k\leq b; \\
\textbf{2)} & b<k\leq c; \\
\textbf{3)} & c<k\leq d. \\
\end{array}%
\right.
\]
In case \textbf{3)} we have $||r_a,r_b,r_c,r_d,r_x|| = n^6\det(A)$
where
\[
\det(A) = \det\smx{
  1 & a & a^2 & a^3 & 0 \\
  1 & b & b^2 & b^3 & 0 \\
  1 & c & c^2 & c^3 & 0 \\
  1 & d & d^2 & d^3 & 1 \\
  1 & x & x^2 & x^3 & \{x \geq k\}  \\
} = \left\{%
\begin{array}{ll}
 \det(A_{5,5})- \det(A_{4,5}) & \text{ if } x\geq k; \\
 - \det(A_{4,5}) & \text{ if } x<k. \\
\end{array}%
\right.
\]
Here $A_{i,j}$ denotes the minor obtained from $A$ by removing the
$i$--th row and the $j$--th column. $A_{4,5}$ and $A_{5,5}$ are
Vandermonde matrices, hence, if $x < k$, $\det(A)$ is equal to the
following polynomial:
\[
\det(A)= -(b-a)(c-a)(c-b)(x-a)(x-b)(x-c),
\]
which is clearly different from zero for $x\neq a,b,c$ and changes
sign whenever $x$  increases and passes through one of the values
$a$, $b$ or $c$, leading to the following signs alternation:
$+\cdots a \cdots - \cdots b \cdots + \cdots c \cdots - \cdots k$.
Then, any lower facet of $\Phi_{k,n}^{(3)}$ should have $b=a+1$
and $c=k-1$; while any upper facet should have $a=0$ and $c=b+1$.
If $x \geq k$, we have the polynomial
\[
\det(A) =
(b-a)(c-a)(c-b)\Big[(d-a)(d-b)(d-c)-(x-a)(x-b)(x-c)\Big],
\]
which is positive whenever $d >x \geq k$, and is negative for
$d<x$, whichever $a$, $b$, $c$ may be. Then, any lower facet
relative to case \textbf{3)} should have $d=n$, and any upper
facet should have $d=k$. So, the facets relative to case
\textbf{3)} are:
\[
\big\{\langle r_t,r_{t+1},r_{k-1},r_n
\rangle\big\}_{t=0,\ldots,k-3}, \qquad \big\{\langle
r_0,r_i,r_{i+1},r_k \rangle\big\}_{i=1,\ldots,k-2}.
\]
The facets relative to  cases \textbf{1)} and \textbf{2)} can be
determined similarly.

\bibliographystyle{plain}
\bibliography{biblio}

\end{document}